\documentclass[10pt,a4paper]{article}
\usepackage{amssymb}
\usepackage[english]{babel}

\newtheorem{theorem}{Theorem}[section]

\newtheorem{lemma}[theorem]{Lemma}
\newtheorem{proposition}[theorem]{Proposition}

\newcommand{\Com}{{\mathbb C}}

\newcommand{\FF}{{\mathbb F}}

\newcommand{\enp}{\begin{flushright} $\Box$ \end{flushright}}

\newcommand{\dis}{\mathbin{\scriptstyle\triangle}}

\begin{document}

\title{From geometry to invertibility preservers%
\thanks{The first author acknowledges the hospitality of the Department of
Mathematics at the University of Ljubljana where most of the work was carried
out. The second author was partially supported by a grant from the Ministry of
Science of Slovenia.}}
\author{Hans Havlicek\\
Institut f\" ur Diskrete Mathematik und Geometrie\\
Technische Universit\" at Wien\\
 Wiedner Hauptstrasse 8--10\\
 A-1040 Wien\\
 Austria\\
havlicek@geometrie.tuwien.ac.at\\[5mm]
{\small and}\\[5mm]
Peter \v Semrl\\
Department of Mathematics\\
        University of Ljubljana\\
        Jadranska 19\\
        SI-1000 Ljubljana\\
        Slovenia\\
peter.semrl@fmf.uni-lj.si}

\date{}
\maketitle

\begin{abstract}
We characterize bijections on matrix spaces (operator algebras) preserving full
rank (invertibility) of differences of matrix (operator) pairs in both
directions.
\end{abstract}
\maketitle



\section{Introduction}

Marcus and Purves \cite{MaP} proved that every unital invertibility preserving
linear map on a matrix algebra is either an inner automorphism or an inner
anti-automorphism. One of the equivalent formulations of Gleason-Kahane-\.
Zelazko theorem \cite{Gle,KaZ,Zel} states that every unital linear functional
on a complex unital Banach algebra ${\cal A}$ sending every invertible element
into a nonzero scalar is multiplicative. Equivalently, if a linear functional
$f:{\cal A} \to \Com$ maps every element $a\in {\cal A}$ into its spectrum
$\sigma (a)$, then $f$ is multiplicative. This two results motivated Kaplansky
to formulate the question under which conditions an invertibility preserving
linear unital map between two algebras must be a Jordan homomorphism
\cite{Kap}. A lot of work has been done on this problem (see the surveys
\cite{Aup,BrS,Sou1}). We will mention here only the results that are relevant
for our paper. Let $X$ be a complex Banach space and $B(X)$ the algebra of all
bounded linear operators on $X$. In 1986 Jafarian and Sourour \cite{JaS} proved
that every surjective unital linear map $\phi : B(X) \to B(X)$ preserving
invertibility in both directions, i.e., having the property that $A$ is
invertible if and only if $\phi (A)$ is invertible, is either of the form $\phi
(A) = TAT^{-1}$, $A\in B(X)$, for some invertible $T\in B(X)$, or of the form
$\phi (A) = TA'T^{-1}$, $A\in B(X)$, for some invertible bounded linear
operator $T: X' \to X$. Here, $A'$ denotes the adjoint of $A$ and $X'$ the dual
of $X$. Under the additional assumption of injectivity the assumption of
preserving invertibility in both directions can be relaxed to the weaker
assumption of preserving invertibility in one direction only \cite{Sou}. The
proof of the result of Jafarian and Sourour was simplified in \cite{Sem}. It is
rather easy to see that a linear map $\phi : B(X) \to B(X)$ is unital and
preserves invertibility in both directions if and only if $\phi$ preserves the
spectrum, that is, $\sigma (\phi (A)) = \sigma (A)$ for every $A\in B(X)$.

An interesting extension of Gleason-Kahane-\. Zelazko theorem was obtained by
Kowalski and S\l odkowski \cite{KoS}. They proved that every functional $f$ on
a complex Banach algebra ${\cal A}$ (they did not assume the linearity of $f$)
satisfying $f(a) - f(b) \in \sigma (a-b)$, $a,b\in {\cal A}$, is linear and
multiplicative up to the constant $f(0)$. Thus, they replaced the two
conditions in Gleason-Kahane-\. Zelazko theorem, the linearity assumption and
the condition $f(a)\in \sigma (a)$, $a\in {\cal A}$, by a single weaker
assumption and got essentially the same conclusion.

In view of this result it is natural to ask if we can do the same with the
above mentioned results on invertibility preserving maps on matrix and operator
algebras. Can we replace the linearity assumption and the invertibility
preserving assumption by a single weaker condition similar to the one in
Kowalski-S\l odkowski theorem? More precisely, can we characterize bijective
maps on matrix algebras and operator algebras satisfying the condition that
$\phi(a) - \phi (b)$ is invertible if and only if $a-b$ is invertible?

The result of Kowalski and S\l odkowski depends heavily on deep results from
analysis. We will answer the above question using the results from geometry. We
should first mention that there is an essential difference between the finite
and the infinite-dimensional case. In the finite-dimensional case our condition
will imply up to a translation the semilinearity of the maps under
consideration, while in the infinite-dimensional case the elementary automatic
continuity methods  will imply the linearity or conjugate-linearity up to a
translation. Moreover, in the finite-dimensional case it makes sense to extend
our result from matrix algebras of square matrices to the spaces of rectangular
matrices. Then, of course, the condition of invertibility will be replaced by
the condition of being of full rank.

Our strategy when considering bijective maps $\phi$ on matrix spaces (operator
algebras) satisfying the condition that $\phi(A) - \phi (B)$ is of full rank
(invertible) if and only if $A-B$ is of full rank (invertible) will be to prove
first that such maps preserve adjacency in both directions. Recall that two
matrices or operators $A$ and $B$ are adjacent if $A-B$ is of rank one. Then we
will apply the so called fundamental theorem of geometry of matrices (or its
analogue for operators) to complete the proof. This connects our results with
the geometry of Grassmann spaces. Let us briefly describe this connection.

Let $M_{m,n}$, $m,n\ge 2$, be the linear space of all $m\times n$ matrices over
a field $\FF$. If $\sigma$ is an automorphism of the field $\FF$ and $A=
[a_{ij}]\in M_{m,n}$ then we denote by $A_\sigma$ the matrix obtained from $A$
by applying $\sigma$ entrywise, $A_\sigma = [\sigma (a_{ij})]$. The fundamental
theorem of geometry of  matrices states that every bijective map $\phi :
M_{m,n}  \to M_{m, n}$ preserving adjacency in both directions is of the form
$A\mapsto TA_\sigma S + R$, where $T$ is an invertible $m\times m$ matrix, $S$
is an invertible $n\times n$ matrix, $R$ is an $m\times n$ matrix, and $\sigma$
is an automorphism of the underlying field. If $m=n$, then we have the
additional possibility that $\phi (A) = TA_{\sigma}^t S + R$ where $T,S,R$ and
$\sigma$ are as above, and $A^t$ denotes the transpose of $A$. This theorem and
its analogues for hermitian matrices, symmetric matrices, and skew-symmetric
matrices were proved by Hua \cite{Hu1}-\cite{Hu8} under some mild technical
assumptions that were later proved to be superfluous (see \cite{Wan}). Let
$m,n$ be integers $\ge 2$. We will consider the Grassmann space whose ``points"
are vector subspaces of $\FF^{m+n}$ of dimension $m$. Chow \cite{Cho} studied
bijective maps on the Grassmann space preserving adjacent pairs of points in
both directions. Recall that $m$-dimensional subspaces $U$ and $V$ are adjacent
if $\dim (U + V) = m+1$. Now, to each $m$-dimensional subspace $U$ of
$\FF^{m+n}$ we can associate an $m\times (m+n)$ matrix whose rows are
coordinates of vectors that form a basis of $U$. Each $m\times (m+n)$ matrix
will be written in the block form $[ X \ Y]$, where $X$ is an $m\times n$
matrix and $Y$ is an $m\times m$ matrix. Two matrices $[ X \ Y]$ and $[ X' \
Y']$ are associated to the same subspace $U$ (their rows represent two bases of
$U$) if and only if $[ X \ Y] = P [ X' \ Y']$ for some invertible $m\times m$
matrix $P$. If this is the case, then $Y$ is invertible if and only if $Y'$ is
invertible. So, we have associated to each point in a Grassmann space a (not
uniquely determined) matrix $[ X \ Y]$. If $Y$ is singular, we say that the
corresponding point in the Grassmann space is at infinity. Otherwise, we
observe that this point can be represented also with the matrix $[Y^{-1}X \ \
\, I]$. The matrix $Y^{-1}X$ is uniquely determined by the point in the
Grassmann space. So, if $U$ and $V$ are two $m$-dimensional subspaces that are
finite points in the Grassmann space, then they can be represented with two
uniquely determined $m\times n$ matrices $T$ and $S$, and it is easy to see
that the subspaces $U$ and $V$ are adjacent if and only if the matrices $T$ and
$S$ are adjacent. Using this connection it is possible to deduce the result of
Chow on bijective maps on a Grassmann space preserving adjacency in both
directions from the fundamental theorem of geometry of matrices (see
\cite{Wan}).

If we consider the special case when $m=n$ and replace in the fundamental
theorem of geometry of matrices the condition of preserving adjacent pairs of
matrices by our assumption of preserving the pairs $A,B$ with the property that
${\rm rank}\, (A-B) =n$, then this corresponds to the study of bijective maps
on the Grassmann space of all vector subspaces of $\FF^{2n}$ of dimension $n$
that preserve the complementarity of subspaces. Such maps were studied by
Blunck and the first author \cite{BlH}. We suspect that this result can be
deduced from our result and the other way around, but we also believe that it
is easier to prove each of them separately. Namely, to prove any of these two
implications seems to be difficult because of the points at infinity.

Now we state our main results. In the finite-dimensional case we will consider
bijective maps on $m\times n$ matrices preserving pairs of matrices whose
difference has a full rank. Of course, if we have such a map $\phi$ then the
map $\psi: M_{n,m} \to M_{n,m}$ defined by $\psi (A) = (\phi (A^t))^t$ has the
same properties. Thus, when studying such maps there is no loss of generality
in assuming that $m\ge n$. We will do this throughout the paper. A matrix $A\in
M_{m,n}$ is said to be of full rank if ${\rm rank}\, A = n$. Let $A,B \in
M_{m,n}$. We write $A\dis B$ if $A-B$ is of full rank.

\begin{theorem}\label{fin}
Let $\FF$ be a field with at least three elements and $m,n$ integers with $m\ge
n \ge 2$. Assume that $\phi : M_{m,n} \to M_{m,n}$ is a bijective map such that
for every pair $A,B \in M_{m,n}$ we have $A\dis  B$ if and only if $\phi (A)
\dis  \phi (B)$. Then there exist an invertible $m\times m$ matrix $T$, an
invertible $n\times n$ matrix $S$, an $m\times n$ matrix $R$, and an
automorphism $\sigma : \FF \to \FF$ such that
$$
\phi (A) = TA_\sigma S + R
$$
for every $A\in M_{m,n}$. If $m=n$, then we have the additional possibility
that $$\phi (A) = TA_{\sigma}^t S + R, \ \ \ A\in M_{n,n},$$ where $T,S,R \in
M_{n,n}$ with $T$ and $S$ invertible, and $\sigma$ is an automorphism of $\FF$.
\end{theorem}

\begin{theorem}\label{infin}
Let $H$ be an infinite-dimensional complex Hilbert space and $B(H)$ the algebra
of all bounded linear operators on $H$. Assume that $\phi : B(H) \to B(H)$ is a
bijective map such that for every pair $A,B \in B(H)$ the operator $A-B$ is
invertible if and only if $\phi (A) - \phi (B)$ is invertible. Then there exist
$R\in B(H)$ and invertible $T,S \in B(H)$ such that either
$$
\phi (A) = TAS + R
$$
for every $A\in B(H)$, or
$$
\phi (A) = TA^t S + R
$$
for every $A\in B(H)$,  or
$$
\phi (A) = TA^* S + R
$$
for every $A\in B(H)$,  or
$$
\phi (A) = T(A^t)^* S + R
$$
for every $A\in B(H)$. Here, $A^t$ and $A^*$ denote the transpose with respect
to an arbitrary but fixed orthonormal basis, and the usual adjoint of $A$ in
the Hilbert space sense, respectively.
\end{theorem}

The converses of both theorems obviously hold true. In the second section we
will prove the finite-dimensional case and in the third one the
infinite-dimensional case. These two sections can be read independently.

\section{The finite-dimensional case}

In this section we will consider matrices over a field $\FF$ with at least
three elements. At a certain point in the proof of our first main theorem we
will identify $m\times n$ matrices with linear operators from $\FF^n$ into
$\FF^m$. For such operators we have the following simple lemma.

\begin{lemma}\label{pez}
Let $T,S : \FF^n \to \FF^m$ be nonzero linear operators and assume that $T$ has
at least two-dimensional image. Then we can find linearly independent vectors
$x,y \in \FF^n$ such that $Tx$ and $Sy$ are linearly independent.
\end{lemma}

{\bf Proof.} Take any $y\in \FF^n$ such that $Sy\not=0$. The set of all vectors
$z\in\FF^n$ with the property that $Tz$ is linearly dependent of $Sy$ is a
proper subspace of $\FF^n$, since the image of $T$ is not contained in the span
of $Sy$. There exist at least two linearly independent vectors of $\FF^n$ which
are not in this particular subspace. One of them is linearly independent of $y$
and gives the required vector $x$. \enp

We have two relations on $M_{m,n}$, that is, the relation of adjacency and
$\dis $. The following result connecting these two relations is the key step in
our proof. We believe it is of some independent interest.

\begin{proposition}\label{wed}
Let $A,B \in M_{m,n}$ be matrices with $A\not=B$. Then the following are
equivalent:
\begin{enumerate}
\item $A$ and $B$ are adjacent.

\item There exists $R\in M_{m,n}$, $R\not=A,B$, such that for every $X\in
M_{m,n}$ the relation $X \dis  R$ yields $X\dis  A$ or $X\dis  B$.
\end{enumerate}
\end{proposition}

{\bf Proof.} Note that none of the above conditions are effected if we replace
$A$ and $B$ by $PAQ-C$ and $PBQ-C$, respectively, where $P$ and $Q$ are
invertible matrices of the appropriate size and $C$ is any $m \times n$ matrix.
Thus if the rank distance between $A$ and $B$ equals $r$ then we may assume
with no loss of generality that $A=0$ and
$$B=\left(\begin{array}{cc}I & 0\\ 0 &0 \end{array} \right)$$
where $I$ is the $r\times r$ identity matrix and the zeros stand for the zero
matrices of the appropriate size.

Assume first that $A$ and $B$ are adjacent. So, without loss of generality, we
have $A=0$ and $B=E_{11}$. Set $R=\lambda E_{11}$, where $\lambda$ is a scalar
different from $0$ and $1$, and $E_{11}$ denotes the matrix with the
$(1,1)$-entry equal to $1$ and all other entries equal to zero. Now let $X\dis
R$. That means that $X-R$ is of full rank or equivalently, the matrix $X-R$
contains at least one invertible $n\times n$ submatrix. We have to consider two
possibilities. Let first assume that one of these submatrices does not contain
the first row. In this case $X$ is of full rank and thus $X\dis A$. Otherwise
any such submatrix contains the first row and we choose one of them. We will
prove that at least one of the corresponding $n\times n$ submatrices of $X-A=X$
and $X-B$ is invertible. So we restrict our attention to these $n\times n$
submatrices. In other words we deal only with the square case $m=n$. Hence
$X-\lambda E_{11}$ is an invertible square matrix. If the first row of
$E_{11}$, i.e.\ $(1,0,\ldots,0)$, is in the subspace spanned by rows
$2,3,\ldots, n$ of $X$ then $X-\lambda E_{11}-\mu E_{11}$ is invertible for all
$\mu\in\FF$, otherwise this holds for all but one $\mu\in\FF$. Therefore
$X-\lambda E_{11} - \mu E_{11}$ is invertible for at least one of the values
$\mu=-\lambda$ or $\mu=-\lambda+1$. Equivalently, at least one of $X=X-A$ or
$X-E_{11}=X-B$ is invertible, as desired. This completes the proof of the first
implication.

To prove the other direction we identify $m\times n$ matrices with linear
operators from $\FF^n$ into $\FF^m$. We assume that $A=0$ and $B: \FF^n \to
\FF^m$ is a linear operator whose image is at least two-dimensional. Let $R:
\FF^n \to \FF^m$ be any linear operator, $R\not=0,B$. We have to find a linear
operator $X: \FF^n \to \FF^m$ such that $X-R$ is injective while $X$ and $X-B$
are not.

The first possibility we will treat is that at least one of the operators $B-R$
and $R$ has rank at least two. Then, by Lemma \ref{pez}, we can find linearly
independent $x,y \in \FF^n$ such that $Bx-Rx$ and $Ry$ are linearly
independent. We first define $X$ on the linear span of $x$ and $y$. We set
$Xx=Bx$ and $Xy=0$. No matter how we will extend $X$ to the whole space these
two equations will guarantee that $X-B$ and $X$ will not be injective. Now,
$(X-R)x = Bx-Rx$ and $(X-R)y = -Ry$ are linearly independent. It is now obvious
that we can extend the linear operator $X$ to the whole space $\FF^n$ such that
the obtained $X-R$ is injective.

It remains to consider the case when both operators $B-R$ and $R$ are of rank
one. By our assumption, $B$ is of rank two. Hence $B= R + (B-R)$ implies that
the ranges of $B-R$ and $R$ meet at $0$ only. So, if we choose any $x,y \in
\FF^n$ such that $(B-R)x\not=0$ and $Ry\not=0$ then $(B-R)x$ and $Ry$ will be
linearly independent. Since $\FF$ has at least three elements, we can choose
these $x$ and $y$ to be linearly independent. Now we can proceed as above. \enp

It is now easy to prove Theorem \ref{fin}. Namely, if $\phi:M_{m,n}\to M_{m,n}$
is a bijective map preserving $\dis$ in both directions, then, by Proposition
\ref{wed} it preserves adjacency in both directions. Thus, the result follows
directly from the fundamental theorem of geometry of matrices.

Observe that in \cite{BlH} there is no need to assume that $\FF$ has at least
three elements, due to the presence of points at infinity.

\section{The infinite-dimensional case}

Let $H$ be an infinite-dimensional complex Hilbert space and $x,y \in H$. The
inner product of $x$ and $y$ will be denoted by $y^* x$. If $x$ and $y$ are
nonzero vectors then $xy^*$ stands for the rank one bounded linear operator
defined by $(xy^*)z = (y^*z)x$, $z\in H$. Note that every bounded rank one
operator can be written in this form. Two operators $A,B \in B(H)$ are said to
be adjacent if $A-B$ is an operator of rank one. We write $A\dis B$ if $A-B$ is
invertible. We start with an analogue of Proposition \ref{wed}.

\begin{proposition}\label{pic}
Let $A,B \in B(H)$ with $A\not= B$. Then the following statements are
equivalent:
\begin{enumerate}
\item $A$ and $B$ are adjacent. \item There exists $R\in B(H)$, $R\not=A,B$,
such that for every $X\in B(H)$ the relation $X \dis R$ yields $X\dis A$ or
$X\dis B$.
\end{enumerate}
\end{proposition}

{\bf Proof.} Note that none of the above conditions are effected if we replace
$A$ and $B$ by $A-C$ and $B-C$, respectively, where $C$ is any bounded linear
operator on $H$. Thus we may assume with no loss of generality that $A=0$.

Assume first that $A=0$ and $B$ are adjacent, that is, $B$ is of rank one. Set
$R=2B$. Suppose that $X-2B$ is invertible. Then
$$
X-2B - \lambda B = ( X- 2B)(I - \lambda (X-2B)^{-1}B)
$$
is invertible if and only if $I - \lambda S$ is invertible, where $S=
(X-2B)^{-1}B$ is an operator of rank one. Every operator of rank one has at
most one non-zero complex number in its spectrum. Hence, $X-2B - (-2B) = X$ is
invertible or $X-2B - (-B) = X-B$ is invertible. This completes the proof of
one direction.

Assume now that $A=0$ and $B$ is an operator whose image is at least
two-dimensional. We have to prove that for every $R\in B(H)$, $R\not= 0,B$,
there exists $X\in B(H)$ such that $X-R$ is invertible and $X$ is singular and
$X-B$ is singular. So, let $R\in B(H) \setminus \{ 0, B \}$.

In the next step we will prove that there exist $x, z\in H$ such that $x$ and
$z$ are linearly independent and $Bz - Rz$ and $Rx$ are linearly independent.
It is enough to show that we can find $x,z\in H$ such that $Bz - Rz$ and $Rx$
are linearly independent. For if $x$ and $z$ are linearly dependent, we can
choose $u\in H$ linearly independent of $x$. Then $z + \lambda u$ and $x$ are
linearly independent for all nonzero $\lambda$ and for all $\lambda$'s small
enough the vectors $B(z + \lambda u ) - R (z + \lambda u) = Bz - Rz + \lambda
(Bu-Ru)$ and $Rx$ are linearly independent as well.

So, let us show that such $x$ and $z$ exist. Assume on the contrary that
$Bz-Rz$ and $Rx$ are linearly dependent for every $x$ and $z$. Then $B-R$ and
$R$ are rank one operators with the same one-dimensional image. It follows that
$B=0$ or $B$ is of rank one, a contradiction.

Now, we define $W$ to be the orthogonal complement of the linear span of $x$
and $z$, where $x$ and $z$ are as in the previous paragraph, and $Z$ to be the
orthogonal complement of $Rx$ and $Bz - Rz$. Then there exists a bounded
invertible linear operator $U: W \to Z$. Define $X\in B(H)$ with
$$
Xx= 0,
$$
$$
Xz = Bz,
$$
and
$$
Xu = Uu + Ru, \ \ \ u\in W.
$$
Because of the first two equations the operators $X$ and $X-B$ are singular.
Since $(X-R) x = -Rx$, $(X-R)z = Bz-Rz$, and $(X-R)u = Uu$, $u\in W$, the
operator $X-R$ is invertible, as desired. \enp

We continue with some technical lemmas.

\begin{lemma}\label{noc}
Let $B,C \in B(H)$. Assume that for every invertible $A\in B(H)$ the operator
$A-B$ is invertible if and only if $A-C$ is invertible. Then $B=C$.
\end{lemma}

{\bf Proof.} Let $\lambda$ be any complex number satisfying
$$
| \lambda | > \| B \| , \| C \|,
$$
and $x,y\in H$ any vectors such that $y^* x =0$. Then $\lambda (I + xy^*)$ is
invertible because $(I + xy^* ) (I - xy^* ) = I$. Hence, $\lambda I + \lambda
xy^* - B$ is invertible if and only if $\lambda I + \lambda xy^* - C$ is
invertible. On the other hand,
$$
\lambda I + \lambda xy^* - B = (I+xy^* )( \lambda I - B +xy^* B) =
(I+xy^*)(I + xy^* B (\lambda I - B)^{-1}) (\lambda I -B)
$$
is invertible if and only if $I + xy^* B (\lambda I - B)^{-1}$ is invertible.
Thus, $I + xy^* B (\lambda I - B)^{-1}$ is invertible if and only if $I + xy^*
C (\lambda I - C)^{-1}$ is invertible, or equivalently, for every scalar
$\lambda$ with $| \lambda | > \| B \| , \| C \|$, and every pair of vectors
$x,y\in H$ with $y^* x = 0$ we have
$$
y^* B(\lambda I -B)^{-1}x = -1 \iff y^* C(\lambda I -C)^{-1}x = -1 .
$$
Fix $\lambda$. Then $y^* T x = 0$ for every pair of orthogonal vectors $x$ and
$y$, where $ T = B(\lambda I -B)^{-1} - C(\lambda I -C)^{-1}$. It follows that
$T = \mu I$ for some scalar $\mu$. Thus, for every $\lambda$ with $| \lambda |
> \| B \| , \| C \|$ we have
$$
B(\lambda I -B)^{-1} - C(\lambda I -C)^{-1} = g(\lambda )I
$$
for some $g(\lambda ) \in \Com$. Obviously, $g(\lambda)$ is holomorphic outside
the circle centered at $0$ with radius $\max \{ \| B \| , \| C \| \}$.
Expressing the above analytic functions with the series and comparing the
coefficients we get
$$
B= C+ \mu I
$$
for some complex number $\mu$. Our assumption implies that
$\sigma(B)\setminus\{0\} = \sigma(C)\setminus\{0\}$. Here $\sigma(B)$ denotes
the spectrum of $B$. It follows that $\mu=0$, as desired. \enp

\begin{lemma}\label{jut}
Let $A,B \in B(H)$ be invertible operators. Assume that for every rank one
operator $xy^* \in B(H)$ the operator $A-xy^*$ is invertible if and only if
$B-xy^*$ is invertible. Then $A=B$.
\end{lemma}

{\bf Proof.} Our assumptions yield that for every pair of vectors $x,y$ the
operator $I - xy^*A^{-1}$ is invertible if and only if $I- xy^*B^{-1}$ is, or
equivalently, $y^* A^{-1}x = 1$ if and only if $y^*B^{-1}x=1$. By linearity we
have $y^* A^{-1}x = y^*B^{-1}x$ for every pair $x,y \in H$, and therefore,
$A^{-1} = B^{-1}$. It follows that $A=B$. \enp

Let us recall that an additive map $T:H\to H$ is called semilinear if there is
an automorphism $\sigma :\Com\to \Com$ such that $T(\lambda x) =
\sigma(\lambda) Tx$ for every $\lambda\in\Com$ and every $x\in H$. Now we are
ready to prove Theorem \ref{infin}.

Let $\phi : B(H) \to B(H)$ be a bijective map such that for every pair $A,B \in
B(H)$ the operator $A-B$ is invertible if and only if $\phi (A) - \phi (B)$ is
invertible. After replacing $\phi$ by $A\mapsto \phi (A) - \phi (0)$ we may
assume that $\phi (0) = 0$. Then $\phi (I)$ is invertible. Replacing $\phi$ by
$A\mapsto \phi(I)^{-1}\phi(A)$ we may further assume that $\phi(I) =I$.

According to Proposition \ref{pic}, $\phi$ preserves adjacency in both
directions. Every rank one operator is adjacent to zero, every rank two
operator is adjacent to some rank one operator, etc. Consequently, $\phi$ maps
the subspace $F(H) \subset B(H)$ of all finite rank operators onto itself. So,
we can apply Theorem 1.5 from \cite{PeS} to conclude that there exist bijective
semilinear maps $T,S : H \to H$ (with the same accompanying automorphism) such
that either $\phi (xy^* ) = (Tx)(Sy)^*$, $x,y \in H$, or $\phi (xy^* ) =
(Sy)(Tx)^*$, $x,y \in H$. The second case can be reduced to the first one if we
replace $\phi$ by $A\mapsto \phi(A)^*$, $A\in B(H)$. So, we may assume that the
first possibility holds true.

Using $\phi (I) = I$ and our assumptions we conclude that $I -xy^*$ is
invertible if and only if $I - (Tx)(Sy)^*$ is invertible, $x,y\in H$. Thus,
$y^* x = 1$ if and only if $(Sy)^* (Tx) = 1$, and by semilinearity,
$$
(Sy)^* (Tx) = 0 \iff y^* x=0, \ \ \ x,y\in H.
$$
Thus, the semilinear maps $T$ and $S$ and their inverses carry closed
hyperplanes (every closed hyperplane is the orthogonal complement of some
nonzero vector) onto closed hyperplanes. Hence, by \cite[Lemmas 2 and 3]{FiL},
$S$ and $T$ are both linear bounded or both conjugate-linear bounded. Thus, we
have $\phi (xy^* ) = T(xy^*)R$, where $T$ and $R=S^*$ are bounded invertible
either both linear, or both conjugate-linear operators. Assume they are both
conjugate-linear. Choosing an orthonormal basis we define $K : H \to H$ to be
the conjugate-linear bijection which maps each vector $x$ into a vector whose
coordinates are obtained from the coordinates of $x$ by complex conjugation. Of
course, $K^2 = I$, the product of two conjugate-linear maps is linear, and
$K(xy^*)K = ((xy^*)^*)^t$, where the transpose is taken with the respect to the
chosen basis. Replacing $\phi$ by $A\mapsto (\phi (A)^t)^*$, $A\in B(H)$, we
reduce the conjugate-linear case to the linear one.

So, we may assume that we have $\phi (xy^* ) = T(xy^*)R$, where $T$ and $R=S^*$
are bounded invertible linear operators. From $(Sy)^* (Tx) =1 \iff y^* x=1$ and
linearity we get $(Sy)^* (Tx) =y^* x$, $x,y\in H$, which further yields that
$T$ is the inverse of $R$. Composing $\phi$ by a similarity transformation we
may further assume that $\phi (xy^* ) = xy^*$, $x,y \in H$.

Let $A\in B(H)$ be invertible.
Applying Lemma \ref{jut} with $B= \phi (A)$ we see that $\phi (A) = A$.

Finally, let $B\in B(H)$ be any operator and set $C=\phi (B)$. Using Lemma
\ref{noc} we conclude that $\phi(B) = B$. This completes the proof. \enp

\end{document}